\documentclass[final]{siamltex}
\usepackage{graphics,graphicx,verbatim}
\usepackage{upquote}
\newlength\myverbindent 
\makeatletter
\def\verbatim@processline{%
 \hspace{\myverbindent}\the\verbatim@line\par}
\makeatother

\def\Im{\hbox{\kern .3pt Im\kern .5pt}}
\def\Re{\hbox{\kern .3pt Re\kern .5pt}}
\def\pput(#1,#2)#3{\noindent\smash{\raise#2pt\hbox to 0pt
   {\kern #1pt #3\hss}}\ignorespaces}

\title{Vandermonde with Arnoldi}

\author{Pablo D. Brubeck, Yuji Nakatsukasa, 
and Lloyd N.~Trefethen\thanks{\texttt{brubeckmarti@maths.ox.ac.uk},
\texttt{nakatsukasa@maths.ox.ac.uk}, and
\texttt{trefethen@maths.ox.ac.uk},
Mathematical Institute, University of Oxford, Oxford, OX2 6GG, UK.}}

\begin{document}

\maketitle

\begin{abstract}
Vandermonde matrices are exponentially ill-conditioned, rendering
the familiar ``polyval(polyfit)'' algorithm for
polynomial interpolation and least-squares fitting ineffective
at higher degrees.
We show that Arnoldi orthogonalization fixes the problem.
\end{abstract}


\pagestyle{myheadings}
\thispagestyle{plain}
\markboth{\sc Brubeck, Nakatsukasa, and Trefethen}
{\sc Vandermonde with Arnoldi}

\section{\label{introd}Introduction}
Fitting polynomials to data by means of Vandermonde matrices is
a notoriously unstable algorithm.  In this note we show how the
problem can be fixed by coupling the Vandermonde construction with
Arnoldi orthogonalization.

Let $x$ be a column vector of $m$ distinct numbers
and $f$ a column vector of $m$ data values.
Given $n\le m-1$, we wish to compute the polynomial
\begin{equation}
p(x) = \sum_{k=0}^n c_k x^n
\end{equation}
with coefficients defined by a system of equations $A c \approx f$,
\begin{equation}
\def\crr{\cr\noalign{\vskip 2pt}}
\def\crrr{\cr\noalign{\vskip 6pt}}
\pmatrix{1 & x_1 & \cdots & x_1^n \crr
1 & x_2 & \cdots & x_2^n \crrr
\vdots & \vdots &  & \vdots \crrr
1 & x_m & \cdots & x_m^n}
\pmatrix{c_0 \crr \vdots \crr c_n}
\approx
\pmatrix{f_1 \crr f_2 \crrr \vdots \crrr f_m}.
\label{system}
\end{equation}
If $n=m-1$, $A$ is square and $\approx$ is interpreted as equality:
$p$ interpolates the data.  If $n<m-1$, $A$ has more rows than
columns, and the system is interpreted in the least-squares sense.
Regardless of $n$, $A$ is of full rank since the points $x_j$
are distinct, so $c$ is uniquely defined.  At least this is true
mathematically, though in cases of interest, $c$ may be nonunique
numerically.

Once (\ref{system}) has been solved, the evaluation of $p$ at a set of
points $s_1,\dots, s_M$ amounts to the computation of the
matrix-vector product
\begin{equation}
\def\crr{\cr\noalign{\vskip 2pt}}
\def\crrr{\cr\noalign{\vskip 6pt}}
\def\crrrr{\cr\noalign{\vskip 10pt}}
\pmatrix{y_1 \crr y_2 \crrrr \vdots \crrrr y_M}
=
\pmatrix{1 & s_1 & \cdots & s_1^n \crr
1 & s_2 & \cdots & s_2^n \crrrr
\vdots & \vdots &  & \vdots \crrrr
1 & s_M & \cdots & s_M^n}
\pmatrix{c_0 \crr \vdots \crr c_n}.
\label{eval}
\end{equation}

\section{\label{ex1}Example 1: interpolation in Chebyshev points}
\begin{figure}
\begin{center}
\vskip .15in
\includegraphics[scale=.84]{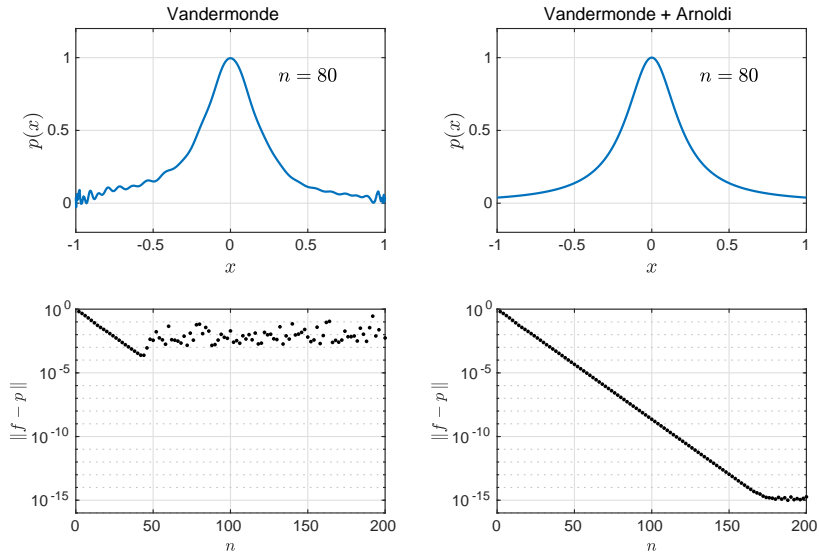}
\end{center}
\caption{\label{fig1}On the left, the degree $n$ 
Chebyshev interpolant to $f(x) = 1/(1+25 x^2)$ computed
unstably by direct application of\/ {\rm (\ref{system})}
and\/ {\rm (\ref{eval})} via the codes
{\tt polyfit} and {\tt polyval} for $n=80$ (above) and its error for even
values of $n$ from $2$ to $200$ (below).
(The results computed by the MATLAB
versions of\/ {\tt polyval} and\/ {\tt polyfit} would be worse.)
On the right, the same computations
with the Arnoldi-based codes {\tt polyfitA} and {\tt polyvalA}.}
\end{figure}

Consider polynomial interpolation of $f(x) = 1/(1+25x^2)$ in $n+1$
Chebyshev points $x_j = \cos(\kern 1pt j\pi/n)$, $0\le j \le n$.
This is a well-conditioned problem, satisfying $\|p \| / \|f\| <
1+(2/\pi) \log (n+1)$ in the $\infty$-norm~\cite[Thm.~15.2]{atap}.
Yet if we apply (\ref{system}) and (\ref{eval}) as written, the
result is a failure.  Specifically, we execute in MATLAB the codes

\begin{verbatim}

function c = polyfit(x,f,n)
A = x.^(0:n);
c = A\f;

\end{verbatim}

\begin{verbatim}
function y = polyval(c,s)
n = length(c)-1;
B = s.^(0:n);
y = B*c;

\end{verbatim}

\noindent where $s$ is a vector of $1000$ points in $[-1,1]$.
The left column of images in Figure~\ref{fig1} shows the results.
For $n=80$, the interpolant is plainly inaccurate, and the plot of
error against $n$ shows that after $n=40$ or so, only 1--3 digits of
accuracy are achieved.  (If one does the same experiment with the
versions of {\tt polyfit} and {\tt polyval} provided in MATLAB,
the results are much worse, a consequence of the order of the
columns of the Vandermonde matrix being reversed.)

On the other hand suppose we execute the following alternative
codes, which will be explained in Section~\ref{VA}.
The ``A'' stands for Arnoldi.

\begin{verbatim}

function [d,H] = polyfitA(x,f,n)
m = length(x);
Q = ones(m,1);
H = zeros(n+1,n);
for k = 1:n
    q = x.*Q(:,k);
    for j = 1:k
        H(j,k) = Q(:,j)'*q/m;
        q = q - H(j,k)*Q(:,j);
    end
    H(k+1,k) = norm(q)/sqrt(m);
    Q = [Q q/H(k+1,k)];
end
d = Q\f;

\end{verbatim}

\begin{verbatim}
function y = polyvalA(d,H,s)
M = length(s);
W = ones(M,1);
n = size(H,2);
for k = 1:n
    w = s.*W(:,k);
    for j = 1:k
        w = w - H(j,k)*W(:,j);
    end
    W = [W w/H(k+1,k)];
end
y = W*d;

\end{verbatim}

\noindent The right column of Figure~\ref{fig1} shows the results,
now converging cleanly down to machine precision.

\section{\label{ex2}Example 2: Least-squares approximation on two intervals}
Our second example moves from interpolation to least-squares.  This
time, the problem is to approximate the function $\hbox{sign}(x)$
on $[-1, -1/3] \cup [1/3,1]$ by a polynomial of degree $n$.
We discretize the domain by 500 equispaced points in each interval
and call {\tt polyfit/polyval} and {\tt polyfitA/polyvalA} to
compute least-squares fits.

\begin{figure}
\begin{center}
\vskip .15in
\includegraphics[scale=.89]{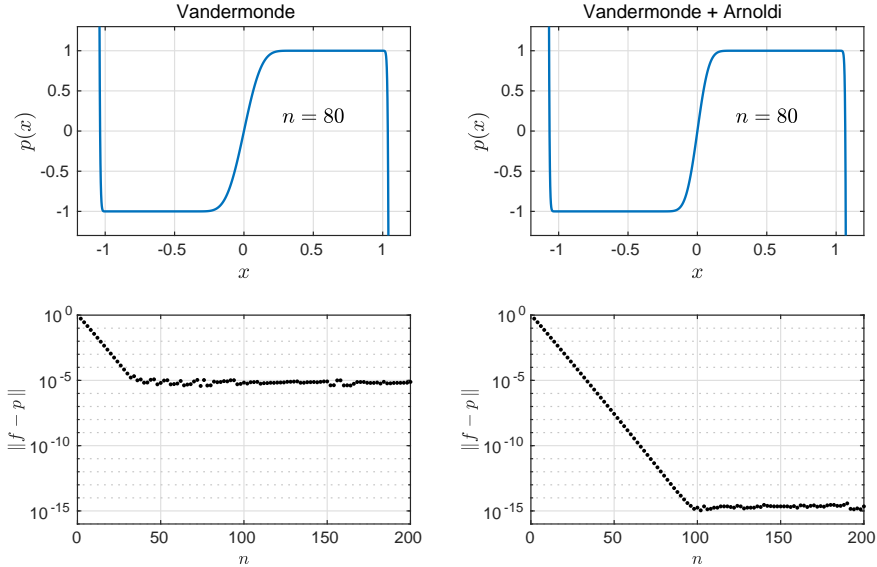}
\end{center}
\caption{\label{fig2}Images as in Fig.~$\ref{fig1}$ but now for
a least-squares problem: polynomial fitting
to $\hbox{sign}(x)$ on $500$ equispaced points
each in the two intervals $[-1, -1/3\kern .5pt ]$ and $[1/3, 1]$.
The unstable algorithm stagnates at\/ $5$ digits of accuracy,
which is enough that to the eye, the computation
appears successful.}
\end{figure}

The results are displayed in Figure~\ref{fig2} in the same
pattern as in Figure~\ref{fig1}, again revealing a great disparity
between algorithms.  A new feature, however, is that the error
with {\tt polyfit/polyval} now goes down to order $O(10^{-5})$
before stagnating at that level.  Ill-conditioned least-squares
computations, at least as carried out by the MATLAB backslash
command, are often sufficiently regularized that the results may
look very good in the ``eyeball norm.''  For some purposes, this
may be all the accuracy one needs, but the figure shows that it
falls far short of what is achievable.

\section{\label{VA}Vandermonde with Arnoldi}
We now explain the Vandermonde with Ar\-noldi algorithm.

First, a word about Vandermonde without Arnoldi, which
solves the equation
\begin{equation}
A c \approx f.
\label{Acf}
\end{equation}
When the code {\tt
polyfit} of Section~\ref{ex1} is executed, exactly what happens
depends on what lies behind the backslash operator---which in the
world of MATLAB is a time-dependent business.  If $A$ is rectangular,
the algorithm will probably be some variant of QR factorization,
with regularizing effects as just mentioned (see \cite{frames} and
also Remark 2.5 of~\cite{barnett}).  If $A$ is square, the algorithm
will probably be Gaussian elimination with partial pivoting.

The introduction of Arnoldi orthogonalization follows a standard idea
in the area of Krylov matrix iterations~\cite{gvl}.  The columns $1,
x, x^2,\dots$ of (\ref{system}) can be regarded as vectors $q_0, Xq_0,
X^2q_0,\dots,$ where $q_0 = (1,\dots , 1)^T$ and $X = \hbox{diag}(x_1,
\dots, x_m)$.  In the Arnoldi process, rather than computing these
vectors as written and then orthogonalizing the result in a QR
factorization, one orthogonalizes at each step to obtain a sequence
of orthogonal vectors $q_0, q_1, q_2,\dots$ spanning the same spaces.
The code {\tt polyfitA} does this in a standard fashion, its only
unusual feature being that $q_k$ is normalized to have
2-norm $m^{1/2}$ rather than $1$, so that the entries individually
are of scale 1 (a matter of convenience).
After $n$ steps, orthogonal vectors $q_0, \dots , q_n$ and an
$(n+1)\times n$ lower-Hessenberg matrix $H$ have been computed such
that
\begin{equation}
X \kern -.3pt Q_- = Q H,
\end{equation}
where $Q$ is the $m\times (n+1)$ matrix with columns $q_0, \dots ,
q_n$ and $Q_-$ is the same matrix without the final column.
Equation (\ref{Acf}) is now solved in its equivalent form
\begin{equation}
Q d\approx f,
\label{Qdf}
\end{equation}
and the resulting
$Q d$ is the $m$-vector of values $p(x_j)$ for the same
polynomial $p$ of degree $n$.  The coefficient
vectors in the two different representations are related mathematically by 
\begin{equation}
d = R^{-1} c
\end{equation}
where $R= Q^T\kern -1.5pt A/m$ is the upper-triangular
matrix such that $A = Q R$.

In the usual Arnoldi application, further operations would now be
carried out on $Q$ and $H$ for purposes such as estimating eigenvalues or
solving systems of equations.
For us, however, all that matters is the
matrix $H$.  The entries of column $k+1$ of $H$ are the coefficients
employed by {\tt polyfitA} to orthogonalize $A q_k$ against $q_0,
\dots, q_k$.  In {\tt polyvalA}, the same operations are applied
with a new operator $S = \hbox{diag}(s_1,\dots , s_M)$ in place of
the original $X = \hbox{diag}(x_1,\dots , x_m)$.
We call the resulting matrices $W_-$ and $W$, satisfying
\begin{equation}
S W_- = W \kern -.9pt H.
\end{equation}
Note that the columns of $W_-$ and $W$ are not orthogonal,
merely (in many cases) approximately so.  If $W\kern -.5pt d$
is now computed, the result is the $M$-vector of values $p(s_j)$.

All this is closely related to orthogonal
polynomials~\cite{davis,gaier,henrici}.  The entries of $H$ are the
recurrence coefficients for a sequence of polynomials orthogonal with
respect to the uniform discrete measure on $\{x_1, \dots, x_m\}$.
If $x$ is a discretization of a real or complex set $\Gamma$,
they can be expected to approximate the recurrence coefficients
for orthogonal polynomials defined by a continuous inner product on
$\Gamma$.  This is what makes the Vandermonde with Arnoldi algorithm
so effective: it is constructing an approximation to orthogonal
polynomials, hence a well-conditioned basis.  The ill-conditioning
of the Vandermonde basis, by contrast, has an elementary explanation
in cases where the points $x_j$ are unequal in size.  If $|x_j|$
varies with $j$, then the powers $|x_j^k|$ vary exponentially.
This means that function information associated with smaller
values of $x_j$ will only be resolvable through exponentially
large expansion coefficients, and accuracy will quickly be lost
in floating point arithmetic.  The wiggles for $|x|\approx 1$ in
the first image of Fig.~\ref{fig1}, for example, are the result of
exponentially large coefficients of terms $x^k$ being required in an
attempt to resolve the function for values $|x|<1$.

Like {\tt polyfit}, {\tt polyfitA} requires $O(mn^2)$ operations,
though with a larger constant.  Subsequent evaluations, however,
cost just $O(Mn)$ operations with {\tt polyval} but $O(Mn^2)$
with {\tt polyvalA}.

Note that since the columns of $Q$ are
orthogonal with norm $m^{1/2}$, the last
line of {\tt polyfitA} could be replaced by \verb|d = Q'*f/m|,
a mathematically equivalent operation involving just $O(mn)$
instead of $O(mn^2)$ operations.  We have not done this since
experiments indicate that accuracy is impaired.  (On the other
hand, this becomes an extremely stable
alternative if each column of $Q$ is computed with not just one but
two rounds of orthogonalization against previous columns.)

If the numbers $x_j$ are real, then $X$ is real symmetric and the
recurrence reduces to three terms at each step: $H$ is tridiagonal
and Arnoldi is equivalent to Lanczos~\cite{gvl}.
In principle this could reduce the operation counts to
$O(mn)$ for {\tt polyfitA} (less than for {\tt polyfit}!) and
$O(Mn)$ for {\tt polyvalA}.\ \ There are probably
applications for which it would be worthwhile to exploit this structure,
but there is a risk of reduced accuracy since the columns of $V$
computed numerically may lose orthogonality.

\section{\label{ex3}Example 3: Fourier extension}
Our third example comes from~\cite{frames}.
Suppose we wish to approximate
$f(x) = 1/(10-9x)$ over the interval $[-1,1]$ using
Fourier series scaled to the larger interval $[-2,2\kern .5pt]$:
\begin{displaymath}
f(x) \approx \sum_{k=-n}^n c_k e^{ik\pi x/2}.
\end{displaymath}
Since $f$ is real, the coefficients will satisfy
$c_{-k} = \overline{c}_k$, and we can write (with
a revised definition of $c_k$)
\begin{displaymath}
f(x) \approx \Re \Bigl(\kern 1pt \sum_{k=0}^n c_k e^{ik\pi x/2}\Bigr).
\label{fourext}
\end{displaymath}
This one-dimensional problem is a model of Fourier extensions in
higher dimensions, which can be used to approximate functions on
irregular domains by embedding them in boxes.  The difficulty is
that the bases are exponentially ill-conditioned.

\begin{figure}
\begin{center}
\vskip .15in
\includegraphics[scale=.89]{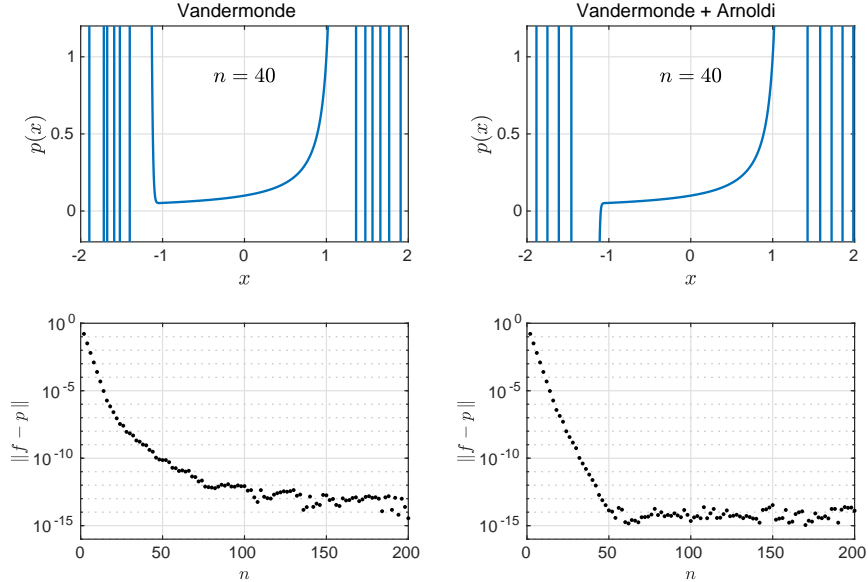}
\end{center}
\caption{\label{fig3}A Fourier extension example
from~{\rm \cite{frames}}, with
$f(x) = 1/(10-9x)$ approximated over $[-1,1]$ by Fourier
series scaled to the larger interval $[-2,2\kern .5pt ]$.  This is equivalent
to approximation by powers $z^k$ over just half of the unit
circle, leading to exponential ill-conditioning of the Vandermonde
matrix.}
\end{figure}

The Vandermonde structure is revealed if we define
$z = e^{i\pi x/2}$, giving
\begin{equation}
f(x) \approx \Re \Bigl(\kern 1pt \sum_{k=0}^n c_k z^k\Bigr), 
\end{equation}
or with $c_k = a_k + i\kern .3pt  b_k$, 
\begin{equation}
f(z) \approx \sum_{k=0}^n \bigl(a_k \Re z^k - b_k \Im z^k \bigr) .
\label{fourext2}
\end{equation}
Our aim is to find
coefficients $\{a_k\}$ and $\{b_k\}$ such that (\ref{fourext2})
is satisfied in the least-squares sense in a set of $m$ points on $\Gamma$
with $m \gg 2n+1$.
This is done by modifying (\ref{system}) to the form
\begin{equation}
\def\crr{\cr\noalign{\vskip 2pt}}
\def\crrr{\cr\noalign{\vskip 6pt}}
\Re \pmatrix{1 & \cdots & z_1^n \crr
1 & \cdots & z_2^n \crrr
\vdots &  & \vdots \crrr
1 & \cdots & z_m^n}
\pmatrix{a_0 \crr \vdots \crr a_n}
-
\Im \pmatrix{z_1 & \cdots & z_1^n \crr
z_2 & \cdots & z_2^n \crrr
\vdots &  & \vdots \crrr
z_m & \cdots & z_m^n}
\pmatrix{b_1 \crr \vdots \crr b_n}
\approx
\pmatrix{f_1 \crr f_2 \crrr \vdots \crrr f_m}.
\label{systemF}
\end{equation}
In {\tt polyfit}, the final line is changed to
\begin{verbatim}

c = [real(A) imag(A(:,2:n+1))]\f;
c = c(1:n+1) - 1i*[0; c(n+2:2*n+1)];

\end{verbatim}
and similarly in {\tt polyfitA} with {\tt d} and {\tt Q} in
place of {\tt c} and {\tt A}.
We also now take the real parts of the results
computed by {\tt polyval} and {\tt polyvalA}.

Figure~\ref{fig3} shows results in the same format as the previous
two examples.  The Fourier extension appears to converge root-exponentially
when implemented with {\tt polyfit} and {\tt polyval}; compare
Fig.~5 of~\cite{frames}.  With {\tt polyfitA} and {\tt
polyvalA}, it converges exponentially.
We note that the success of this computation, and its close relationship
with the previous example after the change of variables $z = e^{i\pi x/2}$,
suggest that in the end there may not be much practical difference 
between Fourier extensions and approximations by ordinary algebraic
polynomials.

\section{\label{ex4}Example 4: Laplace equation and conformal mapping}
Our final example concerns the
solution of the 2D Laplace equation by approximation of the 
solution by the real part of a complex
polynomial~\cite{series}.  To illustrate, let $\Omega$ be the
region bounded by a Jordan curve $\Gamma$ in the complex plane enclosing
the point $0$, and let $g$ be the unique analytic function
that maps $\Omega$ conformally onto the unit disk with $g(0) = 0$ and
$g'(0) > 0$.  We can write 
\begin{equation}
g(z) = z\kern .3pt  e^{\kern .4pt h(z)} = z\kern .3pt  e^{u(z) + iv(z)}
\label{cmap}
\end{equation}
where $h(z) = u(z) + iv(z) $ is an
analytic function with $\Im h(0) = 0$.  Since $|g(z)|= 1$
for $z\in \Gamma$, we have $|h(z)| = -\log(|z|)$ for $z\in\Gamma$, that is,
\begin{equation}
u(z) = -\log(|z|) , \quad z\in \Gamma.
\label{dirichlet}
\end{equation}
Thus the conformal mapping problem reduces to a Laplace problem:
find a harmonic function $u$ in $\Omega$
satisfying the boundary condition~(\ref{dirichlet}),
let $v$ be its harmonic conjugate with $v(0) = 0$, and then $g$ is given
by (\ref{cmap}).
Following~\cite[sec.~2]{conf},
we compute $h$ by approximating it by a polynomial:
\begin{equation}
h(z) \approx \sum_{k=0}^n c_k z^k ,
\end{equation}
or if $c_k = a_k + i\kern .3pt  b_k$, 
\begin{equation}
u(z) \approx \sum_{k=0}^n \bigl(a_k \Re z^k - b_k \Im z^k \bigr) , \quad
v(z) \approx \sum_{k=1}^n \bigl(b_k \Re z^k + a_k \Im z^k \bigr) .
\label{uveq}
\end{equation}
We find coefficients $\{a_k\}$ and $\{b_k\}$ such that (\ref{dirichlet})
is satisfied in the least-squares sense in a set of $m$ points on $\Gamma$
with $m \gg 2n+1$.  This is done again by solving (\ref{systemF}).

\begin{figure}
\begin{center}
\vskip .15in
\includegraphics[scale=.93]{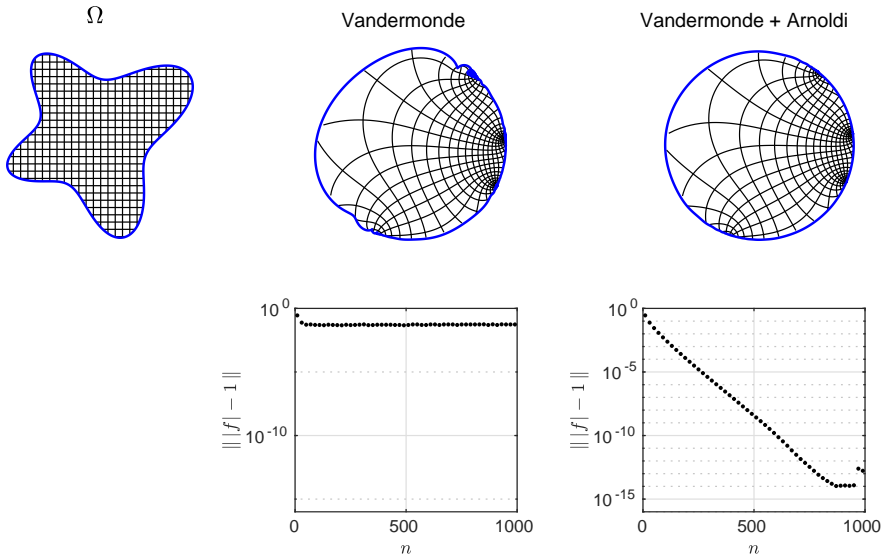}
\end{center}
\caption{\label{fig4}Conformal mapping of a blob onto
the unit disk by the polynomial
expansion method of {\rm (\ref{cmap})--(\ref{uveq})}.  The two
upper-right images correspond to $n=200$.}
\end{figure}

Figure~\ref{fig4} shows the result of such a computation involving
a blob-shaped region $\Omega$.
The Vandermonde solution quickly stagnates, whereas Vandermonde with
Arnoldi is successful.  In this example as in the others, we have
been careful to choose variables $x$ and $z$ with maximum absolute
value equal to~$1$ so that the columns of the matrices (\ref{system})
and (\ref{systemF}) are uniformly scaled.  Without such scaling,
the Vandermonde approach without Arnoldi orthogonalization fails
more dramatically, but that would perhaps be an unfair comparison.

\section{Discussion}
The idea of combining Vandermonde with Arnoldi is not entirely
new.  There is a 1987 paper by Gragg and Reichel~\cite{gragg}
and in several publications, Stylianopoulos and coauthors emphasize the
stability of this approach~\cite{styl,styl13}.  There are also works
related to orthogonal polynomials and quadrature that share
similar mathematics.  Several authors have considered the problem
of constructing orthogonal polynomials for curves or regions in the
complex plane, or families of polynomials defined by criteria
related to orthogonality.  ({\em Bergman polynomials} are defined
by orthogonality with respect to area measure and {\em Szeg\H o
polynomials} by orthogonality over subsets of the unit circle.
Another interesting related recent work is~\cite{faberwalsh}.)
However, the usual emphasis in the literature is on construction
of orthogonal polynomials on a continuous domain as an end
in itself.  Our method instead
constructs discrete orthogonal polynomials on
the fly as a means to enable further computations.  
They are only approximately orthogonal on a continuum,
but in a numerical application, this is all one needs.

There is an interesting literature on Vandermonde matrices, which
are always exponentially ill-conditioned unless the nodes are
uniformly distributed on the unit circle as in the well known case
of the Fast Fourier Transform matrix; for discussion with references,
see~\cite{pan}.  The method presented here bypasses this discussion,
since it bypasses the use of Vandermonde matrices.

For the specific application of polynomial interpolation, we have
not found the Vandermonde with Arnoldi technique in the literature.
Here the use of the barycentric formula is the established stable
technique~\cite{bt,higham}.

Perhaps most important is the rectangular case, the problem of
least-squares fitting with an ill-conditioned basis.  As emphasized
by Adcock and Huybrechs~\cite{frames}, problems of this nature are widespread.
Our own awareness of the importance of Arnoldi orthogonalization
came through the development of ``lightning solvers'' for the Laplace
equation.  Here a polynomial term plays just a supporting role in the
computation, and in writing~\cite{lightning}, we did not recognize
that our computations were sometimes stagnating because this term was
being treated by Vandermonde without Arnoldi.  The {\tt laplace} code
subsequently made available at {\tt people.maths.ox.ac.uk/trefethen/}
remedies the problem by including Arnoldi orthogonalization, as
does the newly introduced
Chebfun code \verb|conformal(...,'poly')|.

Arnoldi stabilization of Vandermonde computations should
have numerous applications beyond those we have illustrated of
interpolation, least-squares fitting, Fourier extension, and solution
of the Laplace equation.  For example, applications
to numerical quadrature could be investigated, as could Arnoldi
stabilization of the calculation of orthogonal polynomials in cases
where that is truly the aim, as well as associated challenges
like the computation of the capacity of a set in the complex plane.
Indeed in any application where a Vandermonde matrix may appear,
it seems likely that introducing Arnoldi orthogonalization may be
a good idea.  It may also be possible to develop
analogous methods for related
applications such as the Laplace eigenvalue problem and the
Helmholtz and biharmonic equations.

\section*{Acknowledgments}
We thank Bernd Beckermann, Lyle Ramshaw, Lothar Reichel, and
Olivier S\`ete for helpful comments.

\indent~~\kern 1in

\end{document}